\documentclass[12pt]{article}
\usepackage{geometry,amsmath,amssymb,bbm,theorem,amscd,graphicx}
\usepackage{setspace} % setspace\UTF{00E3}\UTF{0083}\UTF{0091}\UTF{00E3}\UTF{0083}\UTF{0083}\UTF{00E3}\UTF{0082}¡Þ\UTF{00E3}\UTF{0083}\UTF{0152}\UTF{00E3}\UTF{0082}\UTF{017E}\UTF{00E3}\UTF{0081}\UTF{00AE}\UTF{00E3}\UTF{0082}\UTF{20AC}\UTF{00E3}\UTF{0083}\UTF{00B3}\UTF{00E3}\UTF{0082}\UTF{00AF}\UTF{00E3}\UTF{0083}\UTF{00AB}\UTF{00E3}\UTF{0083}\UTF{0152}\UTF{00E3}\UTF{0083}\UTF{0089}%
\usepackage{colortbl}
\usepackage{ytableau}
\usepackage{tikz}
\usepackage{comment}

\newtheorem{Theorem}{Theorem}[section]

\newtheorem{Example}[Theorem]{Example}

\theoremstyle{remark}
\newtheorem{Remark}{Remark}

\renewcommand\Re{{\mathrm{Re}}}

\newcommand\proof{\it Proof.}

\numberwithin{equation}{section}
\begin{document}
\title{Expressions of content-parametrized Schur multiple zeta-functions via the Giambelli formula\large }
\author{Kohji Matsumoto and Maki Nakasuji
\footnote{The first author is supported by Grants-in-Aid for Scientific Research (B) 
18H01111, and the second author is supported by Grants-in-Aid for Scientific Research (C) 
22K03274.}}
\date{}
\maketitle
%
%
%\pagestyle{empty}
%*********************************************************************
%\address
%{}
%\maketitle
\vskip 1cm
\par\noindent
\begin{abstract}
In this article, we consider the expressions for content-parametrized Schur multiple zeta-functions in terms of multiple zeta-functions of Euler-Zagier type and their star-variants, or in terms of
modified zeta-functions of root systems.
First of all,
we focus on the Schur multiple zeta-function of hook type.
And then, 
applying the Giambelli formula and induction argument,
we obtain the expressions for general content-parametrized Schur multiple zeta-functions.
\end{abstract}

{\small{Keywords: {
Schur multiple zeta-functions, zeta-functions of root systems, 
Euler-Zagier multiple zeta-functions, Giambelli formula}
 

%================Introduction=======================
%%%
%%%
\section{Introduction}\label{sec1}
%%%%%%%%%%%%%%%%%%%%%%%%%%%%%%%%%%%%%%%%%%%%%%%%%%%%%%%%%%%%%%%%%
The multiple zeta-function of Euler-Zagier type and its star-variant are defined by the series
$$\zeta(s_1, \ldots, s_r)=\sum_{0<m_1<\cdots < m_r}\frac{1}{{m_1}^{s_1}\cdots m_r^{s_r}},
\quad
\zeta^{\star}(s_1, \ldots, s_r)=\sum_{0< m_1\leq \cdots \leq  m_r}\frac{1}{{m_1}^{s_1}\cdots m_r^{s_r}},
$$
respectively, where $(s_1, \ldots, s_r)\in {\mathbb C}^r$ with $\Re(s_1), \ldots, \Re(s_{r-1})\geq 1$ and $\Re(s_r)>1$ for convergence.
More precisely, both $\zeta(s_1, \ldots, s_r)$ and $\zeta^{\star}(s_1, \ldots, s_r)$ converge in $\sum_{j=r-i+1}^r{\rm Re}(s_j)>i$ for $1\leq i\leq r$  (see \cite{M}).
Recently, some representation-theoretic generalizations of them have been studied.
One of them is the multiple zeta-functions attached to semisimple Lie algebras, called
the {\it zeta-functions of root systems }introduced by K. Matsumoto, Y. Komori and H. Tsumura 
(\cite{KMT1}). 
It is well known that semisimple Lie algebras over ${\mathbb C}$ are classified by 
root systems coming from their Dynkin diagrams, and
the exact formula for this kind of multiple zeta-function is different for each root system.
For example, the zeta-function of the root system of type $A_r$ is of the following
simple form;
$$
\zeta_r(\underline{\bf s}, A_r)=\sum_{m_1}^{\infty}\cdots \sum_{m_r=1}^{\infty}\prod_{1\leq i<j\leq r+1}(m_i+\cdots + m_{j-1})^{-s(i, j)},
$$
where $s(i, j)$ is the variable corresponding to the root parametrized by $(i, j)$ ($1\leq i, j\leq r+1$, $i\not= j$) and
\begin{align*}
\underline{\bf s}=&(s(1, 2), s(2, 3), \ldots, s(r, r+1), s(1, 3), s(2, 4), \ldots, \\
& s(r-1, r+1), \ldots, s(1, r), s(2, r+1), s(1, r+1)).
\end{align*}
In their papers (\cite{KMT1} and \cite{KMT2}), we can see that various relations among multiple zeta values can be regarded as special cases of functional relations among zeta-functions of root systems. 

The other generalization of multiple zeta-functions of Euler-Zagier type is that associated with combinatorial objects called semi-standard Young tableaux. 
It is called the {\it Schur multiple zeta-functions}, which was introduced by M. Nakasuji, O. Phuksuwan and Y. Yamasaki (\cite{NPY}).
As details are given in the next section, this function has the form of
$$
\zeta_{\lambda}({ \pmb s})=\sum_{M\in \mathrm{SSYT}(\lambda)}\frac{1}{M^{ \pmb s}}, 
$$
where $\lambda$ is a partition and SSYT$(\lambda)$ is a set of semi-standard Young tableaux of shape $\lambda$.
In their paper (\cite{NPY}), after they studied basic properties of this function, 
they obtained some determinant formulas such as Jacobi-Trudi, Giambelli and dual Cauchy formulas under the assumption that $\zeta_{\lambda}({ \pmb s})$
is content-parametrized (in the sense defined in Section \ref{sec2}).
% on variables which we call ``content-parametrized" hereafter.
Furthermore, they investigated skew Schur multiple zeta-functions which is associated with the set difference $\lambda/\mu$ of two partitions $\lambda$ and $\mu$, and quasi-symmetric functions as extensions.

In our previous research(\cite{MN}), the relation between %these two generalizations of multiple zeta-functions are discussed.
zeta-functions of root systems and Schur multiple zeta-functions was discussed.
 There, we obtained various expressions of Schur multiple zeta-functions of anti-hook type 
which are special cases of skew Schur multiple zeta-functions.
One of them is as follows.
\begin{Theorem}\label{firstMN}(\cite[Theorem 3.2]{MN})
For $k, \ell \in {\mathbb N}$, if 
 $\lambda=(
\underbrace{k+1, \cdots, k+1}_{\ell+1 \;{\rm times}}
)$,
  $\mu=
(\underbrace{k, \cdots, k}_{\ell\;{\rm times}}
)$ and
$$
{\pmb s}=
\ytableausetup{boxsize=normal}  
\begin{ytableau}
  \none & \none & \none&   s_{k \ell}\\
  \none & \none & \none&   \vdots\\
  \none & \none & \none &   s_{k1}\\
 s_{00} & s_{10}  & \cdots & s_{k 0}
\end{ytableau}
$$
(note that the way of indexing the variables here is not standard), then
\begin{equation}\label{antihook}
\zeta_{\lambda/\mu}(\pmb s)=\sum_{i=0}^k(-1)^{k-i}\zeta^{\star}(s_{00}, s_{10}, \ldots, s_{i-1, 0}) 
\zeta(s_{k\ell}, s_{k, \ell-1}, \ldots, s_{k0}, s_{k-1, 0}, \ldots, s_{i0})
\end{equation}
holds in the whole space ${\mathbb C}^{k+\ell+1}$, where $\zeta^{\star}=1$ for $i=0$.
\end{Theorem}
Other results in \cite{MN} (Theorem 4.1 in \cite{MN}, for example) are expressions in terms of modified zeta-functions of root systems of type $A$ defined by \eqref{Ahalfstar_def}, \eqref{AH_def} and \eqref{Ahalfstarbullet_def} in Section \ref{sec3}. 

\begin{Remark}\label{WGMDS}
{\rm
One expression among them gives us an analogue of Weyl group multiple Dirichlet series in the sense of Bump, Goldfeld and others (see \cite{B}) and so, this may mean a first link for some undiscovered connections between the theory of Weyl group multiple Dirichlet series and the theory of zeta-functions of root systems which Bump questioned in \cite[p.19]{B}.}
\end{Remark}

Our aim in this article is to obtain the expressions, analogous to the results
proved in \cite{MN}, for Schur multiple zeta-functions of shape $\lambda$.
First of all,
we will focus on the Schur multiple zeta-function of hook type.
In Section \ref{sec3},
we will prove
the expressions similar to \eqref{antihook} in Theorem \ref{firstMN} and those in terms of modified zeta-functions of root systems of type $A$ in \cite{MN}
for hook types.
 % in a much simpler way
 %
 %
%The expressions similar to \eqref{antihook} in Theorem \ref{firstMN} and those in terms of modified zeta-functions of root systems of type $A$ in \cite{MN}
 %can be obtained for the hook types in a much simpler way. %We will discuss it in Section 3.
 %
The aforementioned Giambelli formula for the content-parametrized Schur multiple zeta-functions % under an assumption on variables
is a determinant formula, in which the elements in the matrix in the determinant expression are Schur multiple zeta-functions of hook type.
Therefore, by using this formula, the results in Section \ref{sec3} can be used to
obtain the expressions of
more general content-parametrized Schur multiple zeta-functions.
In Section \ref{sec4}, applying the Giambelli formula with more discussions, we will obtain
the expressions for the content-parametrized Schur multiple zeta-function of shape $\lambda$. 

%In \cite{NPY}, it is known that the Schur multiple zeta-function has a determinant formula called "Giambelli formula".
%In this formula, the elements in the matrix in the determinant expression are hook type Schur multiple zeta-functions, so applying this, we will obtain the expression in terms of multiple zeta- and zet-star functions for the Schur multiple zeta-function of general type.

\vspace{3mm}

{\bf Acknowledgements.} The contents of this article were presented at the ELAZ Conference 2022 in Pozna\'{n}.
The authors are grateful to the organizers of this conference.

%%%%%%%%%%%%%%%%%%%%%%%%%%%%%%%%%%%%%%%%%%%%%%%%%%%%%%%%%%%%%%%%%%%%%%%%%%
%================Section=======================
%
%
\section{Schur multiple zeta-functions and their Giambelli formula}\label{sec2}
%%%%%%%%%%%%%%%%%%%%%%%%%%%%%%%%%%%%%%%%%%%%%%%%
Let $\lambda=(\lambda_1, \ldots, \lambda_m)$ be a partition that is a non-increasing sequence of  a positive integer $n$, i.e. $\lambda_1\geq \lambda_2\geq \cdots \lambda_m>0$ with $\sum_i\lambda_i=n$.
Then a {\it Young diagram} of shape $ \lambda$ is obtained by drawing $\lambda_i$ boxes in the $i$-th row.
The conjugate $\lambda'=(\lambda'_1, \ldots, \lambda'_s)$ of $\lambda$ is defined by $\lambda'_i=\#\{j | \lambda_j\geq i\}$.
In other words, $\lambda'$ is the partition whose Young diagram is the transpose of that of $\lambda$.
Let $T(\lambda,X)$ be the set of all Young tableaux of shape $\lambda$ over a set $X$ and, in particular, $\mathrm{SSYT}(\lambda)\subset T(\lambda,\mathbb{N})$ the set of all semi-standard Young tableaux of shape $\lambda$. Recall that $M=(m_{ij})\in \mathrm{SSYT}(\lambda)$ if and only if $m_{i1}\le m_{i2}\le \cdots$ for all $i$ and $m_{1j}<m_{2j}<\cdots $ for all $j$. For  ${\pmb s}=(s_{ij})\in T(\lambda,\mathbb{C}),$ the Schur multiple zeta-function associated with $\lambda$ is defined as in \cite{NPY} by the series 
$$
\zeta_{\lambda}({ \pmb s})=\sum_{M\in \mathrm{SSYT}(\lambda)}\frac{1}{M^{ \pmb s}}, 
$$
where $M^{ \pmb s}=\displaystyle{\prod_{(i, j)\in \lambda}m_{ij}^{s_{ij}}}$ for $M=(m_{ij})\in \mathrm{SSYT}(\lambda)$ as in Section \ref{sec1}. This series converges absolutely if ${\pmb s}\in W_{\lambda}$ where 
\[
  W_\lambda =
\left\{{\pmb s}=(s_{ij})\in T(\lambda,\mathbb{C})\,\left|\,
\begin{array}{l}
 \text{$\Re(s_{ij})\ge 1$ for all $(i,j)\in \lambda \setminus C(\lambda)$ } \\[3pt]
 \text{$\Re(s_{ij})>1$ for all $(i,j)\in C(\lambda)$}
\end{array}
\right.
\right\}
\]
 with $C(\lambda)$ being the set of all corners of $\lambda$.

%%%%%%%%%%%%%%%%%\section{Giambelli formula}
For a partition $\lambda$, we define two sequences of indices
$p_1, \ldots, p_N$ and $q_1, \ldots, q_N$ by $p_i=\lambda_i-i$ and $q_i=\lambda'_i-i$ for $1\leq i \leq N$
where $N$ is the number of the main diagonal entries of the Young diagram of $\lambda$.
We sometimes write $\lambda=(p_1, \ldots, p_N | q_1, \ldots, q_N)$, which is called 
the Frobenius notation of $\lambda$ (see \cite[Section 1.1]{Mac}).

Let 
$$W_{\lambda}^{\rm{diag}}=\{{\pmb s}=(s_{ij})\in W_{\lambda}\,|\,\text{$s_{ij}=s_{lm}$ if $j-i=m-l$}\}.$$
When ${\pmb s}\in W_{\lambda}^{\rm{diag}}$, we can introduce new variables 
$\{z_k\}_{k\in\mathbb{Z}}$ by the condition $s_{ij}=z_{j-i}$ (for all $i,j$),
and we may regard $\zeta_{\lambda}({\pmb s})$ as a function in variables 
$\{z_k\}_{k\in\mathbb{Z}}$.
We call the Schur multiple zeta-function associated with $\{z_k\}$ {\it content-parametrized Schur multiple zeta-function}, since $j-i$ is named ``content'' (cf. Hamel (\cite{H})).

The following theorem is the Giambelli formula for Schur multiple zeta-functions.

\begin{Theorem}\label{Giambelli} 
(\cite[Theorem 4.5]{NPY})
Let  $\lambda$ be a partition such that $\lambda=(p_1, \cdots , p_N | q_1, \cdots, q_N)$ 
in the Frobenius notation. % where
%$N$ is the number of the elements of the main diagonal.
Assume ${\pmb s}\in W_{\lambda}^{\rm{diag}}$. %has the form $s_{ij}=z_{j-i}$ by a given sequence $\{z_k\}_{k\in {\mathbb Z}}$. 
Then we have
\begin{equation}\label{expij}
\zeta_{\lambda}({\bf s}) = \det(\zeta_{i,j})_{1 \leq i,j \leq N},
\end{equation}
where $\zeta_{i,j}=\zeta_{(p_i+1, 1^{q_j})} ({\bf s}_{ij}^F)$ with ${\bf s}_{ij}^F=
\ytableausetup{boxsize=normal}
  \begin{ytableau}
   z_0 & z_1 & z_2 &\cdots & z_{p_i}\\
   z_{-1}\\
   \vdots \\
    z_{-q_j}
  \end{ytableau}\in W_{(p_i, 1^{q_j})}.
$
%We write $\zeta_{i,0}$ when $\lambda=(p_i+1), {\rm or}\;
%=(p_i-1 | \emptyset )$ in the Frobenius notation.
 \end{Theorem}
 
\begin{Remark}
{\rm The original Giambelli formula is for Schur functions, and the above is its
zeta-analogue. 
In \cite{NPY},
they use the notation $\{a_k\}$ 
 as a variable in $W_{\lambda}^{\rm{diag}}$, 
but we here use the notation $\{z_k\}$, instead.}
\end{Remark}

%%%%%%%%%%%%%%%%%%%%%%%%%%%%%%%%%%%%%%%%%%%%comment %%%%%%%%%%%%
%%%%%%%%%%%%%%%%%%%%     comment   %%%%%%%%%%%comment %%%%%%%%%%%%
%%%%%%%%%%%%%%%%%%%%%%%%%%%%%%%%%%%%%%%%%%%%comment %%%%%%%%%%%%
%%%%%%%%%%%%%%%%%%%%%%%%%%%%%%%%%%%%%%%%%%%%comment %%%%%%%%%%%%
%%%%%%%%%%%%%%%%%%%%%%%%%%%%%%%%%%%%%%%%%%%%comment %%%%%%%%%%%%
\begin{comment}
%\begin{Example}
%When $\lambda=(6,4,4,2,2)$ then $\lambda=(5, 2, 1 | 4, 3, 0)$ is the Frobenius notation, and
%we write 
%$T$
%shortly for $\zeta_{\lambda}(T)$ ($T\in T(\lambda, {\mathbb C})$). Then for ${\pmb s} \in W_{\lambda}^{\rm{diag}}$, we have
%$$
%\zeta_{\lambda}\left(
%\ytableausetup{boxsize=normal}
%  \begin{ytableau}
 %  z_0 & z_1 & z_2 & z_3 & z_4 & z_5 \\
 %  z_{-1} &  z_0  & z_1 & z_2 \\
%   z_{-2} &z_{-1} & z_0& z_1\\
 %   z_{-3 } & z_{-2}\\
 %   z_{-4} & z_{-3}
%  \end{ytableau}
%  \right)
%=
%\det \left|\begin{matrix}
%\zeta_{5,4}
%   & 
%\zeta_{5,3}
%&
%\zeta_{5,0}
%  \\
%\zeta_{2,4}
%  &
%\zeta_{2,3} 
 &
  \zeta_{2,0}
 \\
\zeta_{1,4} 
&
\zeta_{1,3}
&
\zeta_{1,0}
\end{matrix}\right|.$$
\end{Example}
\end{comment}
%%%%%%%%%%%%%%%%%%%%%%%%%%%%%%%%%%%%%%%%%%%%comment %%%%%%%%%%%%
%%%%%%%%%%%%%%%%%%%%%%%%%%%%%%%%%%%%%%%%%%%%comment %%%%%%%%%%%%
%%%%%%%%%%%%%%%%%%%%%%%%%%%%%%%%%%%%%%%%%%%%comment %%%%%%%%%%%%
%%%%%%%%%%%%%%%%%%%%%%%%%%%%%%%%%%%%%%%%%%%%comment %%%%%%%%%%%%

\begin{Example}
{\rm
When $\lambda=(6,4,4,2,2)$ then $\lambda=(5, 2, 1 | 4, 3, 0)$ is the Frobenius notation, and
in this example we write 
$T$
shortly for $\zeta_{\lambda}(T)$ ($T\in T(\lambda, {\mathbb C})$). Then for ${\pmb s} \in W_{\lambda}^{\rm{diag}}$, Theorem \ref{Giambelli} gives}
\vspace{2mm}\\
\hspace{2cm}
\ytableausetup{boxsize=normal}
  \begin{ytableau}
   z_0 & z_1 & z_2 & z_3 & z_4 & z_5 \\
   z_{-1} &  z_0  & z_1 & z_2 \\
   z_{-2} &z_{-1} & z_0& z_1\\
    z_{-3 } & z_{-2}\\
    z_{-4} & z_{-3}
  \end{ytableau}
=
\\
$\det \left|\begin{matrix}

\ytableausetup{boxsize=normal}
  \begin{ytableau}
   z_0 & z_1 & z_2 & z_3 & z_4 & z_5 \\
   z_{-1}  \\
   z_{-2} \\
    z_{-3 }\\
    z_{-4}
  \end{ytableau}

   & 
   \ytableausetup{boxsize=normal}
  \begin{ytableau}
   z_0 & z_1 & z_2 & z_3 & z_4 & z_5 \\
   z_{-1}\\
   z_{-2}\\
   z_{-3}
  \end{ytableau}

&

\ytableausetup{boxsize=normal}
  \begin{ytableau}
 z_0 & z_1 & z_2 & z_3 & z_4 & z_5 
  \end{ytableau}

  \\
  \ytableausetup{boxsize=normal}
  \begin{ytableau}
   z_0 & z_1 & z_2 \\
   z_{-1} \\
   z_{-2} \\
    z_{-3 }\\
    z_{-4}
  \end{ytableau}
  
  &
  
  \ytableausetup{boxsize=normal}
  \begin{ytableau}
   z_0 & z_1 & z_2 \\
   z_{-1}\\
   z_{-2}\\
    z_{-3 }
  \end{ytableau}
  
 &
  
  \ytableausetup{boxsize=normal}
  \begin{ytableau}
   z_0 & z_1 & z_2
  \end{ytableau} 
 
 \\
 
  \ytableausetup{boxsize=normal}
  \begin{ytableau}
   z_0 & z_1\\
   z_{-1}\\
   z_{-2}\\
    z_{-3 }\\
    z_{-4}
  \end{ytableau} 

&

  \ytableausetup{boxsize=normal}
  \begin{ytableau}
   z_0 & z_1\\
   z_{-1}\\
   z_{-2}\\
    z_{-3 }
  \end{ytableau} 

&

  \ytableausetup{boxsize=normal}
  \begin{ytableau}
   z_0 & z_1
  \end{ytableau}

\end{matrix}\right|.$
\end{Example}
%==========================================
%==========================================
%==========================================
%==========================================
%%%%%%%%%%%%%%%%%%%%%%%%%%%%%%%%%%%%%%%%%%%%%%%%%%%%%%%%
\section{Expressions of Schur multiple zeta-functions of hook type}\label{sec3}
%%%%%%%%%%%%%%%%%%%%%%%%%%%%%%%%%%%%%%%%%%%%%%%%%%%%%%%
For $\lambda=(p+1, 1^{q})$, let
\begin{equation}\label{atableau}
{\pmb s}=
%\ytableausetup{boxsize=normal}  
%\begin{ytableau}
%  s_{11} & s_{12} & \cdots & s_{1,p+1}\\
% s_{21} \\
 % \vdots \\
% s_{q+1,1} 
%\end{ytableau}
\begin{array}{|c|c|c|c|}
\hline
s_{11} & s_{12} & \cdots & s_{1,p+1}\\
\hline
 s_{21} \\
\cline{1-1}
\vdots\\
\cline{1-1}
 s_{q+1,1}\\
 \cline{1-1}
\end{array}
=
\ytableausetup{boxsize=normal}  
\begin{ytableau}
\ytableausetup{centertableaux}
  z_{0} & z_{1} & \cdots & z_{p}\\
 z_{-1} \\
  \vdots \\
 z_{-q} 
\end{ytableau}\in W_{\lambda}^{\rm{diag}}\quad, 
\end{equation}
using the notation in the previous section.
The Schur multiple zeta-function of this type is called the {\it Schur multiple zeta-function of hook type}, and we obtain the following expression for it in terms of multiple 
zeta-functions of Euler-Zagier type and their star-variants.

\begin{Theorem}\label{thm3}
For $\lambda=(p+1, 1^{q})$, 
we have
\begin{equation}\label{hook1}
\zeta_{\lambda}({\pmb s})=
\sum_{j=0}^{q}(-1)^j
\zeta^{\star}(z_{-j}, \ldots, z_{-1}, z_0, z_1, \ldots, z_p)
\zeta(z_{-j-1}, \ldots, z_{-q}),
\end{equation}
where we put $\zeta(z_{-j-1}, \ldots, z_{-q})=1$ when $j=q$, and also
\begin{equation}\label{hook2}
\zeta_{\lambda}({\pmb s})=
\sum_{j=0}^{p}(-1)^j
\zeta(z_{j}, \ldots, z_{1}, z_0, z_{-1}, \ldots, z_{-q})
\zeta^{\star}(z_{j+1}, \ldots, z_{p}),
\end{equation}
where we put $\zeta^{\star}(z_{j+1}, \ldots, z_{p})=1$ when $j=p$.
\end{Theorem}

\noindent {\proof}
For the first assertion \eqref{hook1}, 
we carry out the induction for $q$.
In the case $q=0$, it is trivial that the assertion holds.
Assume that the assertion is true for $q-1$, and consider the case of $q$:

\begin{eqnarray}\label{Schurdef}
\zeta_{\lambda}({\pmb s})
&=& \sum_{\substack{m_{11}\leq m_{12}\leq \ldots \leq m_{1,p+1}\\
m_{11}<m_{21}<\ldots <m_{q+1,1}}}
m_{11}^{-z_{0}}m_{12}^{-z_{1}}\ldots m_{1, p+1}^{-z_{p}}
m_{21}^{-z_{-1}}m_{31}^{-z_{-2}}\ldots m_{q+1,1}^{-z_{-q}}.
\end{eqnarray}
The summation on the right-hand side can be divided into two parts:
$$\displaystyle{ \sum_{\substack{m_{11}\leq m_{12}\leq \ldots \leq m_{1, p+1}\\
m_{21}<\ldots <m_{q+1,1}}}
- \sum_{\substack{m_{11}\leq m_{12}\leq \ldots \leq m_{1, p+1}\\
m_{11}\geq m_{21}<\ldots <m_{q+1,1}}}}={\sum}_1 -{\sum}_2, $$
say.
Then obviously
\begin{align*}
 {\sum}_1=\zeta^{\star}(z_0, \ldots, z_p)\zeta(z_{-1}, \ldots, z_{-q}),
\end{align*}
while
\begin{align*}
& {\sum}_2=\sum_{\substack{m_{21}\leq m_{11}\leq m_{12}\leq \ldots \leq m_{1, p+1}\\
m_{21}< m_{31}<\ldots <m_{1, q+1,1}}}=\zeta_{\lambda_- }({\pmb s}^{\flat}),
\end{align*}
where $\lambda_-=(p+2, 1^{q-1})$ and ${\pmb s}^{\flat}=
\begin{array}{|c|c|c|c|c|}
\hline
  s_{21} &  s_{11} & s_{12} & \cdots & s_{1,p+1}\\
\hline
 s_{31} \\
\cline{1-1}
\vdots\\
\cline{1-1}
 s_{q+1,1}\\
 \cline{1-1}
\end{array}
%
%\ytableausetup{boxsize=normal}  
%\begin{ytableau}
%  s_{21} &  s_{11} & s_{12} & \cdots & s_{1, p+1}\\
 %s_{31} \\
 % \vdots \\
 %s_{1, q+1,1} 
%\end{ytableau}
=
\begin{ytableau}
z_{-1} &  z_{0} & z_{1} & \cdots & z_{p}\\
 z_{-2} \\
  \vdots \\
 z_{-q} 
\end{ytableau}
.$

By the assumption for induction, we obtain
\begin{align*}
{\sum}_2 & =\sum_{i=0}^{q-1}(-1)^i \zeta^{\star}(z_{-(i+1)}, \ldots, z_{-2}, z_{-1}, z_{0}, \ldots, z_{p})
\zeta(z_{-(i+2)}, \ldots, z_{-q})\\
& =\sum_{j=1}^{q}(-1)^{j-1} \zeta^{\star}(z_{-j}, \ldots, z_{-2}, z_{-1}, z_{0}, \ldots, z_{p})
\zeta(z_{-(j+1)}, \ldots, z_{-q}).
\end{align*}
This leads to
\begin{align*}
&{\sum}_1-{\sum}_2\\
&=
\zeta^{\star}(z_0, \ldots, z_p)\zeta(z_{-1}, \ldots, z_{-q})
-\sum_{j=1}^{q}(-1)^{j-1} \zeta^{\star}(z_{-j}, \ldots, z_{-2}, z_{-1}, z_{0}, \ldots, z_{p})
\zeta(z_{-(j+1)}, \ldots, z_{-q})\\
&=\sum_{j=0}^{q}
(-1)^{j} \zeta^{\star}(z_{-j}, \ldots, z_{-2}, z_{-1}, z_{0}, \ldots, z_{p})
\zeta(z_{-(j+1)}, \ldots, z_{-q}),
\end{align*}
which completes the proof of \eqref{hook1}. 
The second assertion \eqref{hook2} can be similarly proved by using the induction
for $p$.
$\Box$\vspace{3mm}\\

%%%%%%%%%%%%%%%%%%%%%%%%%%%%%%%%
Theorem \ref{thm3} is a kind of analogue of Theorem \ref{firstMN} in the case of
hook-type.    As mentioned in Section \ref{sec1}, in \cite{MN}, another expression of
$\zeta_{\lambda}({\pmb s})$ in terms of modified zeta-functions of root systems
has been shown.   An analogue of such an expression for the case of hook-type
also exists.

First we have to define the modified zeta-function of root system.
For $r>0$ and $0\leq d \leq r$, we define the modified zeta-function of the root system 
of type $A_r$ by
\begin{align}\label{Ahalfstar_def}
\zeta_{r, d}^{ \bullet}(\underline{\bf s},A_r)=\underbrace{
\left(\sum_{m_1=0}^{\infty}\cdots \sum_{m_d=0}^{\infty}\right)'}_{d \text{ times}}
\underbrace{\sum_{m_{d+1}=1}^{\infty}\cdots \sum_{m_r=1}^{\infty}}_{r-d \text{ times}}\prod_{1\leq i<j\leq r+1}
(m_i+\cdots +m_{j-1})^{-s{(i,j)}},
\end{align}
where 
\begin{multline}\label{ordervector}
\underbar{{\bf s}}=(s{(1,2)},s{(2,3)},\ldots,s{(r,r+1)},s{(1,3)},s{(2,4)},\ldots,s{(r-1,r+1)},\;\ldots,\;
\\
s{(1,r)}, s{(2,r+1)},s{(1,r+1)}).
\end{multline}
with
$s(i,j)$ being the variable corresponding to the root parametrized by
$(i,j)$, and
the prime means that 
the terms $(m_i+\cdots +m_{j-1})^{-s{(i,j)}}$, where $1\leq i<j\leq d+1$ and $m_i=\cdots =m_{j-1}=0$, are omitted.
We also introduce
\begin{align}\label{AH_def}
\zeta_r^H(\underline{\bf s},x, A_r)=\sum_{m_1=1}^{\infty}\cdots \sum_{m_r=1}^{\infty}\prod_{1\leq i<j\leq r+1}
(x+m_i+\cdots +m_{j-1})^{-s{(i,j)}},
\end{align}
and
\begin{align}\label{Ahalfstarbullet_def}
\zeta_{r, d}^{ \bullet, H}(\underline{\bf s}, x, A_r)=\underbrace{
\left(\sum_{m_1=0}^{\infty}\cdots \sum_{m_d=0}^{\infty}\right)}_{d \text{ times}}
\underbrace{\sum_{m_{d+1}=1}^{\infty}\cdots \sum_{m_r=1}^{\infty}}_{r-d \text{ times}}\prod_{1\leq i<j\leq r+1}
(x+m_i+\cdots +m_{j-1})^{-s{(i,j)}},
\end{align}
where $x>0$.
%%%%%%%%%%%%%%%%%%%%%%%%%%%%%%%%5

From the definition of semi-standard Young tableaux, the runnning indices for the Schur multiple zeta-function of hook type  \eqref{Schurdef} satisfy $1\leq m_{11}\leq m_{12}\leq \cdots \leq m_{1, p+1}$ and
$1\leq m_{11}< m_{21}< \cdots < m_{1, q+1,1}$. 
Therefore, setting $m_{12}=m_{11}+a_{1}$ ($a_{1}\geq 0$), 
$m_{13}=m_{11}+a_{1}+a_{2}$ ($a_{1}, a_{2}\geq 0$),
$\cdots$,
\begin{equation}\label{k0ina}
m_{1, p+1}=m_{11}+a_{1}+a_{2}+\cdots +a_{p} \quad (a_{i}\geq 0),
\end{equation}
and  $m_{21}=m_{11}+b_{1}$ ($b_{1}\geq 1$), 
$m_{31}=m_{11}+b_{1}+b_{2}$ ($b_{1}, b_{2}\geq 1$),
$\cdots$,
$$m_{1, q+1,1}=m_{11}+b_{1}+b_{2}+\cdots +b_{q} \quad (b_{j}\geq 1), $$
we can write the Schur multiple zeta-function of hook type as 
\begin{eqnarray}\label{thm4}
\zeta_{\lambda}({\pmb s})
&=& \sum_{\substack{m_{11}\geq 1\\
a_i\geq 0 (1\leq i \leq p)\\
b_j\geq 1(0\leq j\leq q)}}m_{11}^{-z_{0}}(m_{11}+a_1)^{-z_{1}}\cdots (m_{11}+a_1+\cdots + a_{p})^{-z_{p}}
\notag\\
&&\quad\times (m_{11}+b_{1})^{-z_{-1}}
(m_{11}+b_{1}+ b_2)^{-z_{-2}}\cdots 
 (m_{11}+b_{1}+b_{2}+\cdots + b_{q})^{-z_{-q}}
 \notag\\
&=& \sum_{m_{11}\geq 1} m_{11}^{-z_{0}}\zeta^{\bullet, H}_{p, p}({\pmb s}_+, m_{11}, A_p)
 \zeta_{q}^{H}({\pmb s}_-, m_{11}, A_{q}),\label{exp3}
\end{eqnarray}
where ${\pmb s}_+=(z_{1}, z_2, \cdots, z_{p})$ and ${\pmb s}_-=(z_{-1}, z_{-2}, \cdots, z_{-q})$.
We note that this equation is the first step in Theorem \ref{mainthm3} in the next section.

%==========================================
%%%%%%%%%%%%%%%%%%%%%%%%%%%%%%%
%%%%%%%%%%%%%%%%%%%%%%%%%%%%%%%%
%==========================================
\section{Expressions of content-parametrized Schur multiple zeta-functions}\label{sec4}
%%%%%%%%%%%%%%%%%%%%%%%%%%%%%%%%%%%%%%%%%%%%%%%%%%%%%%%%%%%
Now we are ready to prove the main results in the present paper.
Our aim is to generalize the results in Section \ref{sec3} to more general content-parametrized Schur multiple zeta-functions $\zeta_{\lambda}({\pmb s})$ which can be written in terms of $\{z_k\}$ with content $k$. %as in Theorem \ref{Giambelli}.
The results which we will prove in this section are Theorem \ref{mainthm2} and
Theorem \ref{mainthm3}, which are analogues of 
\cite[Theorem 3.2]{MN} and of \cite[Theorem 4.1]{MN}, respectively.

%Let $P=(p_1, p_2, \ldots)$ and $Q=(q_1, q_2, \ldots)$ where $p_i$ and $q_j$ satisfies
%the Frobenius notation $\lambda=(p_1, \cdots , p_N | q_1, \cdots, q_N)$ for a partition $\lambda$.
%Then for $\rho=(1, 2, 3, \ldots)$, we can say $P=\lambda-\rho$ and $Q=\lambda'-\rho$.
First, as a generalization of Theorem \ref{thm3}, we prove the following theorem.
\begin{Theorem}\label{mainthm2}
For the symmetric group ${\frak S}_N$, we have
\begin{align*}
\zeta_{\lambda}({\pmb s})=& \sum_{\sigma\in {\frak S}_N} {\rm{sgn}}(\sigma) \sum_{j_1=0}^{q_1} \cdots \sum_{j_N=0}^{q_N}
(-1)^{j_1+\cdots +j_N}\\
& \times \zeta^{\star}(z_{-j_1}, \ldots, z_0, \ldots, z_{p_{\sigma(1)}})
\zeta^{\star}(z_{-j_2}, \ldots, z_0, \ldots, z_{p_{\sigma(2)}})
\ldots
\zeta^{\star}(z_{-j_N}, \ldots, z_0, \ldots, z_{p_{\sigma(N)}})\\
& \times \zeta(z_{-j_1-1}, \ldots, z_{-q_1})
\zeta(z_{-j_2-1}, \ldots, z_{-q_2})
\ldots
\zeta(z_{-j_N-1}, \ldots, z_{-q_N}).
\end{align*}
\end{Theorem}

\noindent {\proof}
We use the induction for $N$.
When $N=1$, it holds from \eqref{hook1} in Theorem \ref{thm3}.
Assume that the assertion is true for $N-1$, and consider the case of $N$.
The Giambelli formula for the Schur multiple zeta-function \eqref{Giambelli} can be written as
\begin{equation}\label{Delta1}
\zeta_{\lambda}({\pmb s})=\sum_{h=1}^N (-1)^{h+N}\zeta_{h,N}\cdot \Delta_{hN},
\end{equation}
where $\zeta_{i,j}$ is the same notation as in \eqref{expij} and
$$\Delta_{hN}=\det 
\begin{pmatrix}
\zeta_{1,1} & \zeta_{1,2} & \cdots & \zeta_{1,N-1}\\
\vdots & & & \vdots\\
\zeta_{h-1,1} & \zeta_{h-1,2} & \cdots & \zeta_{h-1,N-1}\\
\zeta_{h+1,1} & \zeta_{h+1,2} & \cdots & \zeta_{h+1,N-1}\\
\vdots & & & \vdots\\
\zeta_{N,1} & \zeta_{N,2} & \cdots & \zeta_{N,N-1}\\
\end{pmatrix}
$$
is of size $(N-1)\times (N-1)$. Introducing the notation
$$
\zeta_{i,j}^{\circ}=
\begin{cases}
\zeta_{i,j} & i\leq h-1\\
\zeta_{i+1, j} & i\geq h
\end{cases},
$$
we have
\begin{equation}\label{Deltadet}
\Delta_{hN}=\det 
\begin{pmatrix}
\zeta_{1,1}^{\circ} & \zeta_{1,2}^{\circ} & \cdots & \zeta_{1,N-1}^{\circ}\\
\vdots & & &\vdots \\
\zeta_{N-1,1}^{\circ} & \zeta_{N-1,2}^{\circ} & \cdots & \zeta_{N-1,N-1}^{\circ}\\
\end{pmatrix}.
\end{equation}
Note that $\zeta_{i,j}^{\circ}=\zeta_{i,j}^{\circ}({\pmb s}_{ij}^{\circ})$ where
$${\pmb s}_{ij}^{\circ}={\pmb s}_{i+1,j}^F=
\begin{array}{|c|c|c|c|c|}
\hline
 z_0 & z_1 & z_2 &\cdots & z_{p_{i+1}}\\
 \hline
z_{-1}\\
\cline{1-1}
\vdots\\
\cline{1-1}
 z_{-q_j}\\
 \cline{1-1}
\end{array}
%\ytableausetup{boxsize=normal}
%  \begin{ytableau}
%   z_0 & z_1 & z_2 &\cdots & z_{p_{i+1}}\\
 %  z_{-1}\\
 %  \vdots \\
 %z_{-q_j}
  %\end{ytableau}
  \in W_{(p_{i+1}+1, 1^{q_j})}
$$ 
for $h\leq i\leq N-1$.
We can apply the assumption for the induction to $\Delta_{hN}$
(compare \eqref{expij} and \eqref{Deltadet}).
We have
\begin{align}
\Delta_{hN}&=\sum_{\tau\in \frak{S}_{N-1}}
{\rm{sgn}}(\tau) \sum_{j_1=0}^{q_1} \cdots \sum_{j_{N-1}=0}^{q_{N-1}}
(-1)^{j_1+\cdots +j_{N-1}}\notag\\
& \times \zeta^{\star}(z_{-j_1}, \ldots, z_0, \ldots, z_{p_{\tau'(1)}})
\zeta^{\star}(z_{-j_2}, \ldots, z_0, \ldots, z_{p_{\tau'(2)}})
\ldots
\zeta^{\star}(z_{-j_{N-1}}, \ldots, z_0, \ldots, z_{p_{\tau'({N-1})}})\notag\\
& \times \zeta(z_{-j_1-1}, \ldots, z_{-q_1})
\zeta(z_{-j_2-1}, \ldots, z_{-q_2})
\ldots
\zeta(z_{-j_{N-1}-1}, \ldots, z_{-q_{N-1}}),\label{DeltahN}
\end{align}
where
$$\tau'=\begin{pmatrix}
1 & \cdots & h-1 & h & h+1 & \cdots & N-1 & N\\
1 & \cdots & h-1 & h+1 & h+2 & \cdots & N & h
\end{pmatrix}
\circ
\tau
\in \frak{S}_N.
$$
On the other hand, Theorem \ref{thm3} implies
\begin{align}\label{koremoiru}
\zeta_{h,N}=
\sum_{j_N=0}^{q_N}(-1)^{j_N}
\zeta^{\star}(z_{-{j_N}}, \ldots, z_{-1}, z_0, z_1, \ldots, z_{p_h})
\zeta(z_{-j_N-1}, \ldots, z_{-q_N}).
\end{align}
Substituting \eqref{DeltahN} and \eqref{koremoiru}
into \eqref{Delta1}, we have
\begin{align*}
\zeta_{\lambda}({\pmb s})=& 
\sum_{h=1}^N(-1)^{h+N}
\sum_{j_N=0}^{q_N}(-1)^{j_N}
\zeta^{\star}(z_{-{j_N}}, \ldots, z_{-1}, z_0, z_1, \ldots, z_{p_h})
\zeta(z_{-j_N-1}, \ldots, z_{-q_N})\\
&
\times \sum_{\tau\in {\frak S}_{N-1}} {\rm{sgn}}(\tau) \sum_{j_1=0}^{q_1} \cdots \sum_{j_{N-1}=0}^{q_{N-1}}
(-1)^{j_1+\cdots +j_{N-1}}\\
& \times \zeta^{\star}(z_{-j_1}, \ldots, z_0, \ldots, z_{p_{\tau'(1)}})
\zeta^{\star}(z_{-j_2}, \ldots, z_0, \ldots, z_{p_{\tau'(2)}})
\ldots
\zeta^{\star}(z_{-j_{N-1}}, \ldots, z_0, \ldots, z_{p_{\tau'(N-1)}})\\
& \times \zeta(z_{-j_1-1}, \ldots, z_{-q_1})
\zeta(z_{-j_2-1}, \ldots, z_{-q_2})
\ldots
\zeta(z_{-j_{N-1}-1}, \ldots, z_{-q_{N-1}})\\
=&
\sum_{h=1}^N(-1)^{h+N}
\sum_{\tau\in {\frak S}_{N-1}} {\rm{sgn}}(\tau)
 \sum_{j_1=0}^{q_1} \cdots \sum_{j_{N}=0}^{q_{N}}
(-1)^{j_1+\cdots +j_{N}}\\
& \times \zeta^{\star}(z_{-j_1}, \ldots, z_0, \ldots, z_{p_{\tau'(1)}})
\zeta^{\star}(z_{-j_2}, \ldots, z_0, \ldots, z_{p_{\tau'(2)}})
\ldots
\zeta^{\star}(z_{-j_{N-1}}, \ldots, z_0, \ldots, z_{p_{\tau'(N-1)}})\\
& \times \zeta^{\star}(z_{-j_N}, \ldots, z_0, \ldots, z_{p_h})\\
& \times \zeta(z_{-j_1-1}, \ldots, z_{-q_1})
\zeta(z_{-j_2-1}, \ldots, z_{-q_2})
\ldots
\zeta(z_{-j_{N}-1}, \ldots, z_{-q_{N}}).
\end{align*}
Since  
$\{\tau'\in{\frak S}_{N}\mid \tau\in{\frak S}_{N-1}\}
=\{\sigma\in{\frak S}_{N}\mid \sigma(N)=h\}$,
%$\tau' \in \frak{S}_{N}$ is the permutation with $\tau\in {\frak S}_{N-1}$ and $N \to h$, 
summing over $h$ we have
\begin{align}\label{permutation}
\sum_{h=1}^N(-1)^{h+N}
\sum_{\tau\in {\frak S}_{N-1}} {\rm{sgn}}(\tau) 
=
\sum_{\sigma\in {\frak S}_{N}} {\rm{sgn}}(\sigma). 
\end{align}
This leads to the theorem.
$\Box$

\begin{Remark}
{\rm
If we use \eqref{hook2} in Theorem \ref{thm3} instead of \eqref{hook1}, we obtain an
alternative expression of $\zeta_{\lambda}({\pmb s})$, where the roles of $\zeta$ and
of $\zeta^{\star}$ are reversed.
}
\end{Remark}

Next, as a generalization of \eqref{thm4}, we prove the following theorem, 
which is an expression in terms of modified zeta-functions of root systems.
\begin{Theorem}\label{mainthm3}
For the symmetric group ${\frak S}_N$, we have
\begin{align*}
\zeta_{\lambda}({\pmb s})
=& \sum_{\substack{m_{11}, m_{22}, \ldots, m_{NN}\geq 1}}
(m_{11}\ldots m_{NN})^{-z_{0}}\\
&\times \sum_{\sigma\in {\frak S}_N} {\rm{sgn}}(\sigma)
\prod_{k=1}^N
\zeta^{\bullet, H}_{p_k, p_k}({\pmb s}_{+(k)}, m_{\sigma(k)\sigma(k)}, A_{p_k})
\prod_{j=1}^N \zeta_{q_j}^{H}({\pmb s}_{-(j)}, m_{jj}, A_{q_j}),
\end{align*}
where ${\pmb s}_{+(k)}=(z_{1}, z_2, \cdots, z_{p_k})$ and ${\pmb s}_{-(j)}=(z_{-1}, z_{-2}, \cdots, z_{-q_j})$.
\end{Theorem}

\noindent {\proof}
We use the induction for $N$.
For $N=1$, the assertion is true from \eqref{exp3}.
Assume that it is true for $N-1$ and apply it to $\eqref{Deltadet}$.
Then
\begin{align}
\Delta_{hN}&=
 \sum_{\substack{m_{11}, m_{22}, \ldots, m_{(N-1)(N-1)}\geq 1}}
(m_{11}\ldots m_{(N-1)(N-1)})^{-z_{0}}\notag\\
&\times \sum_{\tau\in {\frak S}_{N-1}} {\rm{sgn}}(\tau)
\prod_{k=1}^{N-1}
\zeta^{\bullet, H}_{p'_k, p'_k}({\pmb s}_{+(k')}, m_{\tau(k)\tau(k)}, A_{p'_k})
\prod_{j=1}^{N-1} \zeta_{q_j}^{H}({\pmb s}_{-(j)}, m_{jj}, A_{q_j}),\label{Dirichlet}
\end{align}
where
$$
k'=\begin{cases}
k & k\leq h-1\\
k+1 & k\geq h
,\end{cases}
\quad\quad
p'_k=\begin{cases}
p_k & k\leq h-1\\
p_{k+1} & k\geq h
.\end{cases}
$$
Dividing the first product in the right-hand side into two parts, we have
\begin{align}
\Delta_{hN}&=
 \sum_{\substack{m_{11}, m_{22}, \ldots, m_{(N-1)(N-1)}\geq 1}}
(m_{11}\ldots m_{(N-1)(N-1)})^{-z_{0}}\notag\\
&\times \sum_{\tau\in {\frak S}_{N-1}} {\rm{sgn}}(\tau)
\prod_{k=1}^{h-1}
\zeta^{\bullet, H}_{p_k, p_k}({\pmb s}_{+(k)}, m_{\tau(k)\tau(k)}, A_{p_k})
\prod_{k=h+1}^{N}
\zeta^{\bullet, H}_{p_k, p_k}({\pmb s}_{+({k})}, m_{\tau(k-1)\tau(k-1)}, A_{p_{k}})\notag\\
&
\times \prod_{j=1}^{N-1} \zeta_{q_j}^{H}({\pmb s}_{-(j)}, m_{jj}, A_{q_j}).\label{DeltahN2}
\end{align}
Substituting \eqref{DeltahN2} and 
$$\zeta_{h,N}=
\sum_{m_{NN}\geq 1}m_{NN}^{-z_0}
\zeta^{\bullet, H}_{p_h, p_h}({\pmb s}_{+(h)}, m_{NN}, A_{p_h})
 \zeta_{q_N}^{H}({\pmb s}_{-(N)}, m_{NN}, A_{q_N})
$$
(which follows from \eqref{thm4})
into \eqref{Delta1}, we have
%\begin{equation}%\label{Delta1}
%\zeta_{\lambda}({\pmb s})=\sum_{h=1}^N (-1)^{h+N}\zeta_{h,N}\cdot \Delta_{hN},
%\end{equation}
%
\begin{align*}
\zeta_{\lambda}({\pmb s})=& 
\sum_{h=1}^N(-1)^{h+N}
 \sum_{\substack{m_{11}, m_{22}, \ldots, m_{NN}\geq 1}}
(m_{11}\ldots m_{NN})^{-z_{0}}
\prod_{j=1}^{N} \zeta_{q_j}^{H}({\pmb s}_{-(j)}, m_{jj}, A_{q_j})
\notag\\
&\times \sum_{\tau\in {\frak S}_{N-1}} {\rm{sgn}}(\tau)
\prod_{k=1}^{h-1}
\zeta^{\bullet, H}_{p_k, p_k}({\pmb s}_{+(k)}, m_{\tau(k)\tau(k)}, A_{p_k})
\zeta^{\bullet, H}_{p_h, p_h}({\pmb s}_{+(h)}, m_{NN}, A_{p_h})\\
& \times \prod_{k=h+1}^{N}
\zeta^{\bullet, H}_{p_k, p_k}({\pmb s}_{+({k})}, m_{\tau(k-1)\tau(k-1)}, A_{p_{k}}).
\end{align*}
Again noting \eqref{permutation}, we obtain the theorem.
%Similar discussion of the last part in the previous theorem leads to the theorem.
%
%
\begin{Remark}
{\rm The right-hand side of \eqref{Dirichlet}  also gives an analogue of Weyl group multiple Dirichlet series in the sense of Bump (\cite{B}); compare with Remark \ref{WGMDS}.}
\end{Remark}

%%%%%%%%%%%%%%%%%%%%%%%%%%%%%%%%5
%
%%%%%%%%%%%%%%%%%%%%%%%%%%%%%%%%%%%%%%%%%%%%%%%%%%%%%
%
%%%%%%%%%%%%%%%%%%%%%%%%%-------------------------
%%%%%%%%%%%%%%%%
%
%
%
%
%==========================================
%==========================================

%
%

%
\bigskip
\noindent
\textsc{Kohji Matsumoto}\\
Graduate School of Mathematics, \\
Nagoya University, \\
Furo-cho, Chikusa-ku, Nagoya, 464-8602, Japan \\
 \texttt{kohjimat@math.nagoya-u.ac.jp}

\medskip

\noindent
\textsc{Maki Nakasuji}\\
Department of Information and Communication Science, Faculty of Science, \\
 Sophia University, \\
 7-1 Kio-cho, Chiyoda-ku, Tokyo, 102-8554, Japan \\
 \texttt{nakasuji@sophia.ac.jp}\\
and\\
Mathematical Institute, \\
Tohoku University, \\
Sendai 980-8578, Japan \\

\begin{thebibliography}{9999999}
%%%%%%%%%%%%%%%%%%%%%%%%%%%%%%%%%%%%%%%%%%%%%%%%%%%%%%%%%%%%%%%%%%%%%%%%%%%%
%%%%%%%%%%%%%%%%%%%%%%%%%%%%%%%%%%%%%%%%%%%%%%%%%%%%%%%%%%%%%%%%%%%%%%%%%%% 
\bibitem[B]{B}
D. Bump, 
Introduction: multiple Dirichlet series, In: Multiple Dirichlet series, $L$-functions and Automorphic Forms, Progress in Mathematics,
{\bf 300}, {\it Burkh\"{a}user/Springer, New York} (2012), 1--36.
%%%%%%%%%%%%%%%%%%%%%%%%%%%%%%%%%%%%%%%%%%%%%%%%%%%%%%%%%%%%%%%%%%%%%%%%%%% 
\bibitem[H]{H}
A. M. Hamel,
Determinantal forms for symplectic and orthogonal Schur functions, 
{\it Canad. J. Math.}, {\bf 49}(2), (1997), 263--282.
%%%%%%%%%%%%%%%%%%%%%%%%%%%%%%%%%%%%%%%%%%%%%%%%%%%%%%%%%%%%%%%%%%%%%%%%%%% 
\bibitem[KMT1]{KMT1}
Y. Komori, K. Matsumoto and H. Tsumura, 
Zeta-functions of root systems,
In: The Conference on $L$-Functions, {\it World Scientific, Hackensack, NJ,} (2007), 115--140.
%%%%%%%%%%%%%%%%%%%%%%%%%%%%%%%%%%%%%%%%%%%%%%%%%%%%%%%%%%%%%%%%%%%%%%%
\bibitem[KMT2]{KMT2}
Y. Komori, K. Matsumoto and H. Tsumura, 
On Witten multiple zeta-functions associated with semisimple Lie algebras. II,
{\it J. Math. Soc. Japan}, {\bf 62}  (2010), 355--394.
%%%%%%%%%%%%%%%%%%%%%%%%%%%%%%%%%%%%%%%%%%%%%%%%%%%%%%%%%%%%%%%%%%%%%%%
\bibitem[Mac]{Mac}
I. G. Macdonald,
Symmetric Functions and Hall Polynomials, 2nd ed.,
Oxford, 1995.
%%%%%%%%%%%%%%%%%%%%%%%%%%%%%%%%%%%%%%%%%%%%%%%%%%%%%%%%%%%%%%%%%%%%%%%%
\bibitem[M]{M}
K. Matsumoto, 
On the analytic continuation of various multiple zeta-functions, In: Number Theory for the Millennium II,
{\it A K Peters Natick, MA}, (2002), 417--440.
%%%%%%%%%%%%%%%%%%%%%%%%%%%%%%%%%%%%%%%%%%%%%%%%%%%%%%%%%%%%%%%%%%%%%%%
\bibitem[MN]{MN}
K. Matsumoto and M. Nakasuji, 
Expressions of Schur multiple zeta-functions of anti-hook type by zeta-functions of root systems,
{\it Publ. Math. Debrecen}, {\bf 98}, (2021), no.3-4, 345--377.
%%%%%%%%%%%%%%%%%%%%%%%%%%%%%%%%%%%%%%%%%%%%%%%%%%%%%%%%%%%%%%%%%%%%%%%
\bibitem[{N}{P}{Y}]{NPY}
M. Nakasuji, O. Phuksuwan, and Y. Yamasaki,
On Schur multiple zeta functions: A combinatoric generalization of multiple zeta functions,
{\it Adv. Math.}, { \bf 333}, (2018), 570--619. 
%%%%%%%%%%%%%%%%%%%%%%%%%%%%%%%%%%%%%%%%%%%%%%%%%%%%%%%%%%%%%%%%%%%%%%%%%%% 

\end{thebibliography}
\end{document}